\documentclass{elsart}

 \usepackage{graphics}
 \usepackage{graphicx}
 \usepackage{epsfig}

\usepackage{amssymb}
\usepackage{amsmath}
\usepackage{amsfonts}

\begin{document}

 \textbf{MANUSCRIPT} \small{published in: Signal Processing 86(10), 2006, pp 2503-3094
  \\ available at: http://dx.doi.org/10.1016/j.sigpro.2006.02.009}  

\begin{frontmatter}

\title{Numerical treatment \\ of an initial-boundary value problem \\
for fractional partial differential equations} 
 \author{Mariusz Ciesielski}
 \ead{mariusz@imi.pcz.pl} \!,
 \author{Jacek Leszczynski}
 \ead{jale@imi.pcz.pl}
 \address{Czestochowa University of Technology\\
Institute of Mathematics and Computer Science\\
ul. Dabrowskiego 73, 42-200 Czestochowa, Poland}

\bigskip

\small{Received 27 April 2005; received in revised form 21 November 2005; \\
accepted 6~Decemeber 2005}

\begin{abstract}
This paper deals with numerical solutions to a~partial differential equation of
fractional order. Generally this type of equation describes a transition from 
anomalous diffusion to transport processes. From a phenomenological point of view,
the equation includes at least two fractional derivatives: spatial and temporal.
In this paper we proposed a new numerical scheme for the spatial derivative,
the so called Riesz-Feller operator. Moreover, using the finite difference
method, we show how to employ this scheme
in the numerical solution of fractional partial differential equations.
In other words, we considered an initial-boundary value problem in one dimensional
space. In the final part of this paper some numerical results and plots 
of simulations are shown as examples. 
\end{abstract}

\end{frontmatter}

\pagestyle{empty}

\section{Introduction}
In the last years fractional derivatives have found numerous 
applications in many fields of physics, mathematics, mechanical engineering,
biology, electrical engineering, control theory and finance 
\cite{Hilf01,Mach01,Main03,Oldh01,Podl01}. One can find interesting properties and
interpretations of fractional calculus in 
\cite{Oldh01,Podl01,Lesz01,Mach02,Mill01,Samk01}. Fractional calculus in mathematics is
a~natural extension of integer-order calculus and gives a useful
mathematical tool for modelling many processes in nature. One of these
processes, in which fractional derivatives have been successfully applied, is
called diffusion. Phenomena which deviate from classical diffusion are described in many papers.
The deviation is called anomalous diffusion. This type of diffusion
is characterized by the nonlinear dependence of the mean square displacement 
$x\left( t\right) $ of a diffusing particle over time~$t$:
$\left\langle \,x^{2}\left( t\right) \right\rangle \sim k_{\gamma }t^{\gamma } $
for $0 < \gamma \leq 2$. This is the opposite of classical diffusion 
where the linear dependence $\left\langle \,x^{2}\left( t\right) 
\right\rangle \sim k_{1}t$ occurs. Analysing changes in the parameter $\gamma$
it may be said that transport phenomena in systems exhibiting sub-diffusion have 
$0<\gamma <1$, and $1<\gamma <2$ in the systems exhibiting super-diffusion.
However the dependence $\left\langle \,x^{2}\left( t\right) 
\right\rangle = \infty $ for $t \rightarrow \infty$ characterizes rare but
extremely large jumps of a diffusing particle - known as L\'{e}vy motion
or L\'{e}vy flights~\cite{Kutn01}. Classical diffusion follows Gaussian statistics and 
Fick's second law for running processes at time $t$, whereas anomalous
diffusion follows non-Gaussian statistics or can be interpreted as the L\'{e}vy
stable densities. The ever increasing amount of literature in which the behaviour of
anomalous diffusion is observed presents some examples of this diffusion in:
disordered vortex lattice for superconductors~\cite{Bouc02}, supercooled
liquids and glasses~\cite{Xia01}, disordered fractal media~\cite{Zeng01},
liquid crystal polymers~\cite{Sant01} and many others.

Using equations with integer-order derivatives to model anomalous diffusion 
in the above processes may not reflect their real real behaviour
Therefore anomalous diffusion is described
by a spatio-temporal fractional partial differential equation in which classical spatial and
temporal derivatives are replaced by derivatives of fractional order.
The anomalous diffusion equation includes the classical diffusion equation, 
and therefore, this equation has a general form.
Moreover, the anomalous diffusion equation may also describe wave propagation or 
advection processes.

Equations of anomalous diffusion with time and/or space fractional derivatives
have been proposed and analysed by numerous authors, for example Nigmatullin in 
\cite{Nigm01,Nigm02}, Bouchaud in \cite{Bouc01}, Wyss and Schneider in
\cite{Schn01,Schn02,Wyss01}, West in \cite{West01}, 
Gorenflo, Fabritiis and Mainardi~in \cite{Gore02},
Gorenflo and Mainardi in~\cite{Gore04}, Mainardi in~\cite{Main02},
Metzler and Klafter~\cite{Metz01}, Hilfer in \cite{Hilf01} and recently by Agrawal in 
\cite{Agra01}. Nevertheless, the theoretical analysis and numerical methods 
applied to solve fractional diffusion equations present difficulties.
In many papers, the autors have considered and solved problems in the infinite domain.
Hilfer \cite{Hilf01} and Klafter and Metzler \cite{Klaf01} described the analytical 
solution to these equations in terms of Fox's $H$-function. 
In \cite{Agra01} Agrawal presented an analytical solution over time for
the anomalous diffusion equation with boundary conditions of 
the first kind. He based his approach on the Laplace transform in terms of the Mittag-Leffer
function. However, the numerical approximation for the series of expansions of
these functions are a little problematic, especially for greater values of
function arguments.

Some numerical methods such as the finite difference method (FDM) and
the finite element method (FEM) are more suitable for solving the anomalous diffusion equation
in more general (non-linear) form. The numerous works by 
Gorenflo and Mainardi \cite{Gore02,Gore04}, Podlubny \cite{Podl01} and many others should be noted. 
The difference scheme for fractional derivatives is based on the definition in the
Gr\H{u}nwald-Letnikov form~\cite{Podl01,Samk01}. This from can make the scheme more flexible and straightforward.
It is also a difficult task to solve the boundary value problem of these equations.
Ciesielski and Leszczynski in \cite{Cies01}, Yuste in \cite{Yust01}
proposed and analysed different cases of numerical solutions to the fractional
diffusion equation where specific kinds of boundary conditions were selected.

In this paper, we present a solution over space to a fractional partial differential equation
in which the Riesz-Feller potential is included. In several papers~\cite{Deng01,Meer01} are
shown some numerical schemes which didn't analyse originally proposed this potential.
The authors introduced their own potentials which are suitable from a phenomenological
point of view.
Numerical treatments of space fractional differential equations are not as popular
in literature as the numerical solution to time fractional differential equations.
This arises because the discretization of the
Riesz-Feller potential is more difficult due to the appearance of singularities in
this operator. On the basis of the Gr\H{u}nwald-Letnikov definition of fractional
derivative, Gorenflo and Mainardi \cite{Gore02,Gore04} proposed a~method of 
discretization for this operator in the infinite domain. However, this 
descretization does not provide continuity for all values of the derivative order.
Taking into account this inconvenience we shall propose another discrete
form of the Riesz-Feller operator in order to restrict its continuity
for arbitrary values of the derivative order. 
We use a new approach to solve numerically an~equation of anomalous diffusion.

\section{Mathematical model}
In this paper we consider an equation in the following form 
\begin{equation}
\frac{\partial }{\partial t}C(x,t)=K_{\alpha }\frac{\partial ^{\alpha }}
{\partial \left\vert x\right\vert ^{\alpha } _{\theta }}C(x,t) \text{, } \quad  \quad 
t\geq 0 \text{, }  \quad  x\in \mathbb{R} \text{, }
\label{e_ad}
\end{equation}
where 
$C(x,t)$ is a field variable, $\left( \partial ^{\alpha }
/ \partial \left\vert x\right\vert ^{\alpha } _{\theta } \right) C(x,t)$ is the Riesz-Feller fractional
operator~\cite{Chec01,Metz01,Samk01}, 
$\alpha $ is the real order of this operator, $K_{\alpha }$ is the
coefficient of generalized (anomalous) diffusion with the unit of measurement 
$ m^{\alpha }/s $.
According to~\cite{Gore02,Gore04} the Riesz-Feller fractional operator 
for $0<\alpha \leq 2$, $\alpha \neq 1$, and for the one-variable function $u(x)$ is
defined as
\begin{equation}
\frac{\partial ^{\alpha }}{\partial \left\vert x\right\vert ^{\alpha } _{\theta }} u(x) \equiv 
 {}_{x}D_{\theta }^{\alpha }u\left( x\right) =
 -\left[ c_{L}\left( \alpha ,\theta \right) \,_{-\infty }D_{x}^{\alpha }u\left( x\right) 
+c_{R}\left( \alpha ,\theta \right) \,_{x}D_{+\infty}^{\alpha } u\left( x\right) \right] 
\text{,}
\label{e_RF}
\end{equation}
where 
\begin{equation}
{}_{-\infty }D_{x}^{\alpha }u\left( x\right) =
\left( \dfrac{d}{dx} \right)^m \left[ _{-\infty }I_{x}^{m-\alpha }u\left( x\right) \right] ,
\label{e_DL}
\end{equation}
\begin{equation}
{}_{x}D_{+\infty }^{\alpha }u\left( x\right) =%
\left( -1\right)^m \left( \dfrac{d}{dx} \right)^m 
 \left[ _{x}I_{+\infty }^{m-\alpha }u\left( x\right) \right] ,
\label{e_DR}
\end{equation}
with $m\in \mathbb{N}$, $m-1<\alpha \leq m$ as the left-side and right-side 
Riemann-Liouville fractional derivatives.
In formula (\ref{e_RF}) coefficients $c_{L}\left( \alpha ,\theta \right) $, $c_{R}\left( \alpha
,\theta \right) $ (for $0<\alpha \leq 2$, $\alpha \neq 1$, and for $%
\left\vert \theta \right\vert \leq \min \left( \alpha ,2-\alpha \right) $)
have the following forms
\begin{equation}
c_{L}\left( \alpha ,\theta \right) =\frac{\sin \dfrac{\left( \alpha -\theta
\right) \pi }{2}}{\sin \left( \alpha \pi \right) }\text{, }\quad c_{R}\left(
\alpha ,\theta \right) =\frac{\sin \dfrac{\left( \alpha +\theta \right) \pi 
}{2}}{\sin \left( \alpha \pi \right) }\text{.}  
\label{e_C}
\end{equation}

In expressions (\ref{e_DL}) and (\ref{e_DR}) the fractional operators 
$_{-\infty}I_{x}^{\alpha }u\left( x\right) $
and $_{x}I_{\infty }^{\alpha }u\left(x\right) $
are defined as the left- and right-side of Weyl fractional integrals
~\cite{Gore02,Gore04,Oldh01,Podl01,Samk01}. Thus we have
\begin{align}
_{-\infty }I_{x}^{\alpha }u\left( x\right) & = \frac{1}{\Gamma (\alpha )}
\int_{-\infty }^{x} \frac{u\left( \xi \right) }{(x-\xi )^{1-\alpha }}d\xi 
\text{,} \label{e_IL} \\
_{x}I_{\infty }^{\alpha }u\left( x\right) & = \frac{1}{\Gamma(\alpha )}
\int_{x}^{\infty } \frac{u\left( \xi \right) }{(\xi-x)^{1-\alpha }}d\xi 
\text{.} \label{e_IR}
\end{align}

Assuming $\alpha =2$ in~(\ref{e_ad}) we obtain classical diffusion,
i.e. the heat transfer equation. On the other hand, the classical transport equation is obtained if $\alpha =1$
and the parameter of skewness $\theta $ in~(\ref{e_C}) admits extreme values. 
Taking into consideration the above changes in parameter $\alpha $
we assume its variations within the range $0<\alpha \leq 2$. 
Analysing in~(\ref{e_ad}) the behaviour of parameter $\alpha$ $(\alpha < 2)$
we observed some transition between the transport and propagation processes. 

We used the Green's function~\cite{Gore04} in the analytical solution of Eq.~(\ref{e_ad}).
However, Eq.~(\ref{e_ad}) needs  a numerical solution
when an additional non-linear term may occur. Some numerical methods used 
to solve fractional partial differential equations can be found in~\cite{Gore02}.
However, these methods apply the infinite domain without boundary conditions. 

In this work, we will consider Eq.~(\ref{e_ad}) in the one dimensional domain 
$\Omega :L\leq x\leq R$ with boundary-value conditions of the first kind 
(Dirichlet conditions) as
\begin{equation}
\left\{ 
\begin{array}{ll}
x=L: \quad & C\left( L,t\right) =g_{L}\left( t\right) \\ 
x=R: \quad & C\left( R,t\right) =g_{R}\left( t\right) 
\end{array}
\right. \text{ for } t>0
\label{e_boundary}
\end{equation}
and with the initial-value condition
\begin{equation}
\left. C\left( x,t\right) \right\vert _{t=0}=c_{0}\left( x\right) \text{.}
\label{e_init}
\end{equation}


\section{Numerical method}
According to the FDM \cite{Ames01,Hoff01} we consider a~discrete form of~Eq.~(\ref{e_ad})
both in time and space. In our previous work \cite{Cies01} we solved numerically 
an anomalous diffusion equation similar to~(\ref{e_ad}) 
where only the time-fractional derivative was taken into account.
We called this method the Fractional FDM (FFDM).
Extending our considerations, we can say that the solution to equation~(\ref{e_ad})
needs proper approximation of the Riesz-Feller derivative~(\ref{e_RF}) in numerical schemes. 

Here, we introduce a definition of the fractional derivative in the
Caputo form \cite{Capu01,Podl01} as
\begin{equation}
{}_{a}^{C}D_{x}^{\alpha }u\left( x\right)  =
{}_{a}I_{x}^{m-\alpha }u^{\left( m \right) }\left( x\right) , 
\label{e_DLC} 
\end{equation}
\begin{equation}
{}_{x}^{C}D_{b}^{\alpha }u\left( x\right)  =
{}_{x}I_{b}^{m-\alpha }u^{\left( m \right) }\left( x\right) , 
\label{e_DRC}
\end{equation}
where $m\in \mathbb{N}$, $m-1<\alpha \leq m$, $a,b \in \mathbb{R}$, $a < x$, $b > x$ 
and $u^{\left( m \right) }$ 
are first and second derivatives, for $m = 1,2$. 
The following relations are present between the Riemann-Louville and the Caputo forms:
\begin{eqnarray}
{}_{a}D_{x}^{\alpha }u\left( x\right) &=& 
 {}_{a}^{C}D_{x}^{\alpha }u\left( x\right) + 
 \sum_{k=0}^{m-1}\left. \frac{d^{k}}{dx^{k}}u \left( x\right)
\right\vert _{x = a}\frac{(x-a)^{k-\alpha  }}{\Gamma \left( k-\alpha  +1\right) }%
\text{,}  \label{eq_rel_Cap_Riem_L}
\\
{}_{x}D_{b} ^{\alpha }u\left( x\right)  &=&
 {}_{x}^{C}D_{b}^{\alpha }u\left( x\right)  +
 \sum_{k=0}^{m-1}\left. \frac{d^{k}}{dx^{k}}u \left( x\right)
\right\vert _{x = b}\frac{(b-x)^{k-\alpha  }}{\Gamma \left( k-\alpha  +1\right) }%
\text{,}  \label{eq_rel_Cap_Riem_R}
\end{eqnarray}
Based on the assumptions $\left\vert \lim_{a\rightarrow -\infty }\left. \left(
d^{k}/dx^{k}\right) u \left( x\right) \right\vert _{x=a}\right\vert <\infty $ and as well \linebreak
$\left\vert \lim_{b\rightarrow \infty }\left. \left(
d^{k}/dx^{k}\right) u \left( x\right) \right\vert _{x=b}\right\vert <\infty $ for $k=0,...,m-1$,
the terms occuring on the left sides of (\ref{eq_rel_Cap_Riem_L}) and (\ref{eq_rel_Cap_Riem_R})
tend to zero. Thus, when the lower/upper limit of integration 
tends to minus/plus infinity we have
\begin{equation}
_{-\infty }D_{x}^{\alpha }u\left( x\right) ={}_{-\infty }^{\hspace{0.35cm}C}D_{x}^{\alpha
}u\left( x\right) \quad \text{and}\quad _{x}D_{+\infty }^{\alpha }u\left( x\right)
={}_{x}^{C}D_{+\infty }^{\alpha }u\left( x\right) \text{.}  
\label{e_LIU}
\end{equation} 
 
\subsection{Approximation of the Riesz-Feller derivative}
As we start numerical analysis from the discretization of 
operators~(\ref{e_DLC}) and~(\ref{e_DRC}) respectively,
we introduce a homogenous spatial grid $-\infty <\ldots
<x_{i-2}<x_{i-1}<x_{i}<x_{i+1}<x_{i+2}<\ldots <\infty $ with the step 
$h=x_{k}-x_{k-1}$. We denote the value of function $u\left(x\right)$ 
at the point $x_{k}$ as $u_{k}=u\left( x_{k}\right) $,
for $k\in \mathbb{Z}$. We take into account only the function of one variable 
in order to simplify notations and we denote 
$c_{L}=c_{L}\left( \alpha ,\theta \right) $, 
$c_{R}=c_{R}\left( \alpha ,\theta \right) $.
In accordance with changes in parameter $\alpha$ in~(\ref{e_ad}) we distinguish
two cases of discrete approximation of the Riesz-Feller derivative.

The first case includes changes in parameter $\alpha$ for 
the range $0<\protect\alpha <1$.  
We rewrite operator~(\ref{e_RF}) using the Caputo definition as 
\mathindent0cm
\begin{eqnarray}
&&_{x_{i}}^{C}D_{\theta }^{\alpha }u\left( x_{i}\right)  =-\left[ c_{L}\,\frac{1%
}{\Gamma (1-\alpha )}\int\limits_{-\infty }^{x_{i}}\frac{u^{\prime }\left(
\xi \right) }{(x_{i}-\xi )^{\alpha }}d\xi -{}{}c_{R}\frac{1}{\Gamma
(1-\alpha )}\int\limits_{x_{i}}^{\infty }\frac{u^{\prime }\left( \xi \right) 
}{(\xi -x_{i})^{\alpha }}d\xi \,\right]   \notag \\
&& \quad =\frac{-1}{\Gamma (1-\alpha )}\left[ c_{L}\sum\limits_{k=0}^{\infty
}\int\limits_{x_{i-k-1}}^{x_{i-k}}\frac{u^{\prime }\left( \xi \right) }{%
(x_{i}-\xi )^{\alpha }}d\xi \,-{}{}c_{R}\sum\limits_{k=0}^{\infty
}\int\limits_{x_{i+k}}^{x_{i+k+1}}\frac{u^{\prime }\left( \xi \right) }{(\xi
-x_{i})^{\alpha }}d\xi \,\right] \label{e_n01a} \\ 
&& \quad \approx \frac{-1}{\Gamma (1-\alpha )}\left[ c_{L}\sum\limits_{k=0}^{\infty
}\widetilde{u}_{i-k}^{\prime }\int\limits_{x_{i-k-1}}^{x_{i-k}}\frac{1}{%
(x_{i}-\xi )^{\alpha }}d\xi \,-{}{}c_{R}\sum\limits_{k=0}^{\infty }%
\widetilde{\widetilde{u}}_{i+k}^{\prime }\int\limits_{x_{i+k}}^{x_{i+k+1}}%
\frac{1}{(\xi -x_{i})^{\alpha }}d\xi \,\right] 
\notag 
\end{eqnarray}%
where $\widetilde{u}_{j}^{\prime }$ and $\widetilde{\widetilde{u}}_{j}^{\prime }$
are difference schemes which approximate the first
derivative of integer order in the intervals $\left[ x_{j-1},x_{j}\right] $ and 
$\left[ x_{j},x_{j+1}\right] $, respectively. We propose the following
weighed forms of these difference schemes as
\mathindent2cm
\begin{eqnarray}
\widetilde{u}_{j}^{\prime } &=&\frac{1}{2}\left( \frac{u_{j}-u_{j-1}}{h}+
\frac{\left( 1-\lambda _{1}\right) \left( u_{j}-u_{j-1}\right) +\lambda
_{1}\left( u_{j+1}-u_{j}\right) }{h}\right)   \notag \\
&=&\frac{1}{2h}\left[ \lambda _{1}u_{j+1}+2\left( 1-\lambda _{1}\right)
u_{j}+\left( \lambda _{1}-2\right) u_{j-1}\right] ,
\label{e_n01_1L}
\end{eqnarray}%
\begin{eqnarray}
{\widetilde{\widetilde{u}}}_{j}^{\prime } &=&\frac{1}{2}\left( \frac{%
u_{j+1}-u_{j}}{h}+\frac{\left( 1-\lambda _{1}\right) \left(
u_{j+1}-u_{j}\right) +\lambda _{1}\left( u_{j}-u_{j-1}\right) }{h}\right)  
\notag \\
&=&\frac{1}{2h}\left[ \left( 2-\lambda _{1}\right) u_{j+1}+2\left( \lambda
_{1}-1\right) u_{j}+\left( -\lambda _{1}\right) u_{j-1}\right] ,
\label{e_n01_1R}
\end{eqnarray}%
where $\lambda _{1}=\lambda _{1}\left( \alpha ,\theta \right)
=\alpha -\left\vert \theta \right\vert $, $\lambda _{1} \in \left[ 0,1\right] $. 
We introduce these formulae because we want to obtain various transitions between the difference schemes
which are connected with the first derivative
of integer order. For example, after putting $\lambda _{1}=1$ into~(\ref{e_n01_1L}) 
and~(\ref{e_n01_1R}) we obtain wide known the central-difference approximation 
of first derivative, and after putting $\lambda _{1}=0$ we get 
the backward- (\ref{e_n01_1L}) or forward- (\ref{e_n01_1R}) difference schemes.
In this way we would like to avoid the problem of singularity in the operator
(\ref{e_RF}) for $\alpha \rightarrow 1^{-}$ and $\theta = 0$
(in further calculations of the discrete form of the Riesz-Feller operator we obtain
finite values of the coefficents for this case). Simultaneously, 
we can obtain classical schemes used for hyperbolic equations
of first order for $\alpha \rightarrow 1^{-}$ and $\theta \rightarrow 1^{\pm }$,
e.g. $ \left( u_{i+1} - u_{i} \right) / h$ and $ \left( u_{i} - u_{i-1} \right) / h$.

Denoting by
\mathindent1cm
\begin{eqnarray}
v_{k}^{\left( \alpha \right) } &=&\frac{1}{\Gamma (1-\alpha )}%
\int\limits_{x_{i-k-1}}^{x_{i-k}}\frac{1}{(x_{i}-\xi )^{\alpha }}d\xi \,=%
\frac{1}{\Gamma (1-\alpha )}\int\limits_{x_{i+k}}^{x_{i+k+1}}\frac{1}{(\xi
-x_{i})^{\alpha }}d\xi \,  \notag \\
&=&h^{1-\alpha }\frac{\left( k+1\right) ^{1-\alpha }-k^{1-\alpha }}{\Gamma
\left( 2-\alpha \right) },
\label{e_n01_v}
\end{eqnarray}%
we have
\begin{equation}
\begin{split}
_{x_{i}}D_{\theta }^{\alpha }u_{i}\approx & -\left[ c_{L}\sum\limits_{k=0}^{%
\infty }\,\frac{1}{2h}\left[ \lambda _{1}u_{i-k+1}+2\left( 1-\lambda
_{1}\right) u_{i-k}+\left( \lambda _{1}-2\right) u_{i-k-1}\right]
v_{k}^{\left( \alpha \right) }\right.  \\
& \left. -c_{R}\sum\limits_{k=0}^{\infty }\frac{1}{2h}\left[ \left(
2-\lambda _{1}\right) u_{i+k+1}+2\left( \lambda _{1}-1\right) u_{i+k}+\left(
-\lambda _{1}\right) u_{i+k-1}\right] v_{k}^{\left( \alpha \right) }\right] .
\end{split}
\label{e_n01b}
\end{equation}

Finally we can write the discrete form of (\ref{e_RF}) as 
\begin{equation}
_{x_{i}}D_{\theta }^{\alpha }u_{i}\approx \dfrac{1}{h^{\alpha }}
\sum\limits_{k=-\infty }^{\infty }u_{i+k} w_{k} \text{,}  \label{e_n01_scheme}
\end{equation}
where coefficients $w_{k} = w_{k} \left( \alpha , \theta \right) $ have the following form:
\mathindent1cm
\begin{eqnarray}
&&w_{k} =\dfrac{-1}{2\Gamma \left( 2-\alpha \right) 
}   
\label{e_n01_w} \\
&&\times \left\{ 
\begin{array}{ll}
\left( \left( \left\vert k\right\vert +2\right) ^{1-\alpha }\lambda
_{1}+\left( \left\vert k\right\vert +1\right) ^{1-\alpha }\left( 2-3\lambda
_{1}\right) \right.  &  \\ 
\left. \quad \quad +\left\vert k\right\vert ^{1-\alpha }\left( 3\lambda
_{1}-4\right) +\left( \left\vert k\right\vert -1\right) ^{1-\alpha }\left(
2-\lambda _{1}\right) \right) c_{L} & \text{ for }k\leq -2, \\ 
\left( 3^{1-\alpha }\lambda _{1}+2^{1-\alpha }\left( 2-3\lambda _{1}\right)
+3\lambda _{1}-4\right) c_{L}+\lambda _{1}c_{R} & \text{ for }k=-1, \\ 
\left( 2^{1-\alpha }\lambda _{1}-3\lambda _{1}+2\right) \left(
c_{L}+c_{R}\right) & \text{ for }k=0, \\ 
\left( 3^{1-\alpha }\lambda _{1}+2^{1-\alpha }\left( 2-3\lambda _{1}\right)
+3\lambda _{1}-4\right) c_{R}+\lambda _{1}c_{L} & \text{ for }k=1, \\ 
\left( \left( k+2\right) ^{1-\alpha }\lambda _{1}+\left( k+1\right)
^{1-\alpha }\left( 2-3\lambda _{1}\right) \right.  &  \\ 
\left. \quad \quad +k^{1-\alpha }\left( 3\lambda _{1}-4\right) +\left(
k-1\right) ^{1-\alpha }\left( 2-\lambda _{1}\right) \right) c_{R} 
& \text{ for }k\geq 2.
\end{array}
\right.  \notag
\end{eqnarray}

The second case involves changes in parameter $\alpha$ 
for the range $1<\protect\alpha \leq 2$. As in previous calculations,
we write operator~(\ref{e_RF}) using the Caputo form as 
\mathindent0cm
\begin{eqnarray}
_{x_{i}}^{C}D_{\theta }^{\alpha }u\left( x_{i}\right)  &=&- \! \left[ c_{L}
\frac{1}{\Gamma (2-\alpha )}\int\limits_{-\infty }^{x_{i}}
\frac{u^{\prime \prime }\left( \xi \right) }{(x_{i}-\xi )^{\alpha -1}}d\xi 
+ c_{R}\frac{1}{\Gamma (2-\alpha )}\int\limits_{x_{i}}^{\infty }
\frac{u^{\prime \prime }\left( \xi \right) }{(\xi -x_{i})^{\alpha -1}}d\xi \right] 
 \notag \\
&\approx &\frac{-1}{\Gamma (2-\alpha )} 
\left[ c_{L}\sum\limits_{k=0}^{\infty }\widetilde{u}_{i-k}^{\prime \prime }
\int\limits_{x_{i-k-1}}^{x_{i-k}}\frac{1}{(x_{i}-\xi )^{\alpha -1}}d\xi \right. \notag \\
&& \hspace{2cm} \left. +\, c_{R}\sum\limits_{k=0}^{\infty } \widetilde{\widetilde{u}}_{i+k}^{\prime \prime }
\int\limits_{x_{i+k}}^{x_{i+k+1}}\frac{1}{(\xi -x_{i})^{\alpha -1}}d\xi 
\right] 
\label{e_n12a}
\end{eqnarray}%
where $\widetilde{u}_{j}^{\prime \prime }$ and $\widetilde{\widetilde{u}}%
_{j}^{\prime \prime }~$ are difference schemes of the second derivative
of integer order which we approximate by the following formulae
\begin{eqnarray}
\widetilde{u}_{j}^{\prime \prime } &=&\frac{1}{2} \!
\left( \frac{u_{j+1}\! -\! 2u_{j}\! + \! u_{j-1}} {h^{2}}
+\frac{\left( 1\! - \! \lambda _{2}\right) \!
\left( u_{j+1}\! -\! 2u_{j}\! + \! u_{j-1}\right) 
+\lambda _{2}\! \left( u_{j}\! - \! 2u_{j-1} \! + \! u_{j-2}\right) }{h^{2}}\! \right)   \notag \\
&=&\frac{1}{2h^{2}}\left[ \left( 2-\lambda _{2}\right) u_{j+1}+\left(
3\lambda _{2}-4\right) u_{j}+\left( 2-3\lambda _{2}\right) u_{j-1}+\lambda
_{2}u_{j-2}\right] ,
\label{e_n12_2L}
\end{eqnarray}%
\begin{eqnarray}
\widetilde{\widetilde{u}}_{j}^{\prime \prime } &=&\frac{1}{2} \! 
\left( \frac{u_{j+1}\! -\! 2u_{j}\! + \! u_{j-1}}{h^{2}}
+\frac{\left( 1\! -\! \lambda _{2}\right) \!
\left( u_{j+1}\! -\! 2u_{j}\! + \! u_{j-1}\right) +\lambda _{2}\! 
\left( u_{j+2}\! -\! 2u_{j+1}\! + \! u_{j}\right) }{h^{2}} \! \right)   \notag \\
&=&\frac{1}{2h^{2}}\left[ \lambda _{2}u_{j+2}+\left( 2-3\lambda _{2}\right)
u_{j+1}+\left( 3\lambda _{2}-4\right) u_{j}+\left( 2-\lambda _{2}\right)
u_{j-1}\right] ,
\label{e_n12_2R}
\end{eqnarray}%
where $\lambda _{2}=\lambda _{2}\left( \alpha ,\theta \right)
=2-\left( \alpha +\left\vert \theta \right\vert \right) $, $\lambda _{2} \in \left[ 0,1\right] $. 
By putting $\lambda _{2}=0$ into~(\ref{e_n12_2L}) and~(\ref{e_n12_2R}) we
determine the classical central-difference scheme, and for $\lambda _{2}=1$
we obtain the backward/forward four-point discretizations of the second derivative
of integer order.

Denoting
\begin{eqnarray}
v_{k}^{\left( \alpha \right) } &=&\frac{1}{\Gamma (2-\alpha )}%
\int\limits_{x_{i-k-1}}^{x_{i-k}}\frac{1}{(x_{i}-\xi )^{\alpha -1}}d\xi \,=%
\frac{1}{\Gamma (2-\alpha )}\int\limits_{x_{i+k}}^{x_{i+k+1}}\frac{1}{(\xi
-x_{i})^{\alpha -1}}d\xi \,  \notag \\
&=&h^{2-\alpha }\frac{\left( k+1\right) ^{2-\alpha }-k^{2-\alpha }}{\Gamma
\left( 3-\alpha \right) },
\label{e_n12_v}
\end{eqnarray}%
we have
\begin{eqnarray}
_{x_{i}}D_{\theta }^{\alpha }u_{i} &\approx &-\left[ c_{L}\sum%
\limits_{k=0}^{\infty }\,\frac{1}{2h^{2}}\left[ \left( 2-\lambda _{2}\right)
u_{i-k+1}+\left( 3\lambda _{2}-4\right) u_{i-k}+\left( 2-3\lambda
_{2}\right) u_{i-k-1}\right. \right.   \notag \\
&&\left. +\lambda _{2}u_{i-k-2}\right] v_{k}^{\left( \alpha \right)
}+c_{R}\sum\limits_{k=0}^{\infty }\frac{1}{2h^{2}}\left[ \lambda
_{2}u_{i+k+2}+\left( 2-3\lambda _{2}\right) u_{i+k+1}\right.   \notag \\
&&\left. \left. +\left( 3\lambda _{2}-4\right) u_{i+k}+\left( 2-\lambda
_{2}\right) u_{i+k-1}\right] v_{k}^{\left( \alpha \right) }\right] .
\label{e_n12b}
\end{eqnarray}

Finally we can write the discrete form of expression~(\ref{e_RF}) as 
\begin{equation}
_{x_{i}}D_{\theta }^{\alpha }u_{i}\approx \dfrac{1}{h^{\alpha }}%
\sum\limits_{k=-\infty }^{\infty }u_{i+k} w_{k} 
\text{,}
\label{e_n12_scheme}
\end{equation}%
where coefficients $w_{k} = w_{k} \left( \alpha , \theta \right) $ are 
\mathindent0cm
\begin{eqnarray}
&&w_{k} = \dfrac{-1}{2\Gamma \left( 3-\alpha \right) 
}
\label{e_n12_w} \\
&&\times \left\{ 
\begin{array}{ll}
\left( \left( \left\vert k\right\vert +2\right) ^{2-\alpha }\left( 2-\lambda
_{2}\right) +\left( \left\vert k\right\vert +1\right) ^{2-\alpha }\left(
4\lambda _{2}-6\right) +\left\vert k\right\vert ^{2-\alpha }\left(
6-6\lambda _{2}\right) \right.  &  \\ 
\left. \quad \quad +\left( \left\vert k\right\vert -1\right) ^{2-\alpha
}\left( 4\lambda _{2}-2\right) +\left( \left\vert k\right\vert -2\right)
^{2-\alpha }\left( -\lambda _{2}\right) \right) c_{L} & \hspace{-0.8cm} \text{for }k\leq -2,
\\ 
\left( 3^{2-\alpha }\left( 2-\lambda _{2}\right) +2^{2-\alpha }\left(
4\lambda _{2}-6\right) -6\lambda _{2}+6\right) c_{L}+\left( 2-\lambda
_{2}\right) c_{R} & \hspace{-0.8cm} \text{for }k=-1, \\ 
\left( 2^{2-\alpha }\left( 2-\lambda _{2}\right) +4\lambda _{2}-6\right)
\left( c_{L}+c_{R}\right) & \hspace{-0.8cm} \text{for }k=0,\\ 
\left( 3^{2-\alpha }\left( 2-\lambda _{2}\right) +2^{2-\alpha }\left(
4\lambda _{2}-6\right) -6\lambda _{2}+6\right) c_{R}+\left( 2-\lambda
_{2}\right) c_{L} & \hspace{-0.8cm} \text{for }k=1, \\ 
\left( \left( k+2\right) ^{2-\alpha }\left( 2-\lambda _{2}\right) +\left(
k+1\right) ^{2-\alpha }\left( 4\lambda _{2}-6\right) +k^{2-\alpha }\left(
6-6\lambda _{2}\right) \right.  &  \\ 
\left. \quad \quad +\left( k-1\right) ^{2-\alpha }\left( 4\lambda
_{2}-2\right) +\left( k-2\right) ^{2-\alpha }\left( -\lambda _{2}\right)
\right) c_{R} & \hspace{-0.8cm} \text{for }k\geq 2.
\end{array}
\right.  \notag
\end{eqnarray}

Summarising the above calculations, we presented difference schemes for 
the Riesz-Feller fractional derivative.
It should be noted that expressions~(\ref{e_n01_scheme}) 
and~(\ref{e_n12_scheme}) are represented by the weighted sum over 
discrete values of function $u$ at all the node's points.
If index $k$ tends to $0$, i.e. in the nearest proximity of 
an arbitrary point $x_i$, one may observe higher values of $w_{k}$.
Whereas values $w_{k}$ decrease to zero 
for those nodes furthest from the point $x_i$. Table~\ref{tab1} shows sample values 
of $w_{k}$ calculated for different 
values of parameter $\alpha$. Here we assumed $\theta =0$, $c_L = c_R$, 
and some symmetry in coefficients $w_{k} = w_{|k|}$, for $k \neq 0$.
It should be noted that coefficients $w_{k}( 2,0)$
are identical as to those for the central difference approximation for the
second derivative, for $\alpha = 2$.
Assuming the skewness parameter to be $\theta = 1^{\pm}$ and $\alpha = 1^{\pm}$ 
we can see that coefficients $w_{k}(1^{\pm}, 1^{\pm})$ tend
to values represented by the backward/forward difference scheme 
for the first derivative.

\begin{table}[h]
\caption{Coefficients $w_{k} \left( \alpha ,0 \right) $ being dependent on
parameter $\alpha $ and $k$ respectively}
\begin{tabular}{|l|c|c|c|c|c|}
\hline
$k$ & \multicolumn{5}{|c|}{$\alpha $} \\ \cline{2-6}
& $0.1$ & $0.5$ & $1^{-}$ & $1.5$ & $2$ \\ \hline
\multicolumn{1}{|r|}{$0$} & \multicolumn{1}{|r|}{$-0.993029$} & 
\multicolumn{1}{|r|}{$-0.963132$} & \multicolumn{1}{|r|}{$-0.857606$} & 
\multicolumn{1}{|r|}{$-1.498970$} & \multicolumn{1}{|r|}{$-2$} \\ \hline
\multicolumn{1}{|r|}{$1$} & \multicolumn{1}{|r|}{$0.041819$} & 
\multicolumn{1}{|r|}{$0.170296$} & \multicolumn{1}{|r|}{$0.253710$} & 
\multicolumn{1}{|r|}{$0.574964$} & \multicolumn{1}{|r|}{$1$} \\ \hline
\multicolumn{1}{|r|}{$2$} & \multicolumn{1}{|r|}{$0.022853$} & 
\multicolumn{1}{|r|}{$0.067624$} & \multicolumn{1}{|r|}{$0.064577$} & 
\multicolumn{1}{|r|}{$0.125442$} & \multicolumn{1}{|r|}{$0$} \\ \hline
\multicolumn{1}{|r|}{$3$} & \multicolumn{1}{|r|}{$0.014264$} & 
\multicolumn{1}{|r|}{$0.036213$} & \multicolumn{1}{|r|}{$0.029047$} & 
\multicolumn{1}{|r|}{$0.020048$} & \multicolumn{1}{|r|}{$0$} \\ \hline
\multicolumn{1}{|r|}{$4$} & \multicolumn{1}{|r|}{$0.010322$} & 
\multicolumn{1}{|r|}{$0.023595$} & \multicolumn{1}{|r|}{$0.016789$} & 
\multicolumn{1}{|r|}{$0.009118$} & \multicolumn{1}{|r|}{$0$} \\ \hline
\multicolumn{1}{|r|}{$5$} & \multicolumn{1}{|r|}{$0.008054$} & 
\multicolumn{1}{|r|}{$0.016974$} & \multicolumn{1}{|r|}{$0.010996$} & 
\multicolumn{1}{|r|}{$0.005125$} & \multicolumn{1}{|r|}{$0$} \\ \hline
\multicolumn{1}{|r|}{$10$} & \multicolumn{1}{|r|}{$0.003751$} & 
\multicolumn{1}{|r|}{$0.006116$} & \multicolumn{1}{|r|}{$0.002926$} & 
\multicolumn{1}{|r|}{$0.000906$} & \multicolumn{1}{|r|}{$0$} \\ \hline
\end{tabular}
\label{tab1}
\end{table}


\subsection{Complete fractional FDM}
Although discretization of the Riesz-Feller derivative in space was proposed,
we describe the FDM for the equation of anomalous diffusion~(\ref{e_ad})
in this subsection. Here, we restrict the numerical
solution to only one-dimensional space. In comparison to the classical
diffusion equation where discretization of the second derivative over space is
approximated by the central difference scheme, the anomalous 
diffusion equation requires generalized schemes given by 
formulae~(\ref{e_n01_scheme}) and~(\ref{e_n12_scheme}), respectively.
Boundary conditions have a direct influence to the numerical solution 
not only on boundary nodes but also in internal nodes of the domain. 

We introduce a temporal grid $0=t^{0}<t^{1}<\ldots <t^{f}<t^{f+1}<\ldots < \infty $
with the grid step $\Delta t=t^{f+1}-t^{f}$. 
At a point $x_{k}$ at the moment of time $t^{f}$ we denote
the function $C\left( x,t\right) $ as $C_{k}^{f}=C\left( x_{k}, t^{f}\right) $,
for $k\in \mathbb{Z}$ and $f\in \mathbb{N}$. 


\subsubsection{Pure initial value problem}
In the explicit scheme of FDM we replaced~(\ref{e_ad}) with the following formula 
\begin{equation}
\dfrac{C_{i}^{f+1}-C_{i}^{f}}{\Delta t}=K_{\alpha }\dfrac{1}{h^{\alpha }}%
\sum\limits_{k=-\infty }^{\infty }C_{i+k}^{f}w_{k} 
\text{.}
\label{e_FDMp1}
\end{equation}%
After simplifications we obtained the final form as
\begin{equation}
C_{i}^{f+1}=\sum\limits_{k=-\infty }^{\infty }C_{i+k}^{f}p_{k} \text{,}
\label{e_FDMp2}
\end{equation}%
where coefficients $p_{k} = p_{k} \left( \alpha , \theta  \right) $ are 
\begin{equation}
p_{k} =\left\{ 
\begin{array}{ll}
1+K_{\alpha }\dfrac{\Delta t}{h^{\alpha }}w_{0} & 
\text{ for }k=0\text{,}  \\ 
K_{\alpha }\dfrac{\Delta t}{h^{\alpha }}w_{k} & 
\text{ for }k\neq 0 \text{.}
\end{array}
\right.
\label{e_FDM_p}
\end{equation}

Now we calculate the sum of all coefficent values $p_{k}$:
\begin{equation}
\sum_{k=-\infty }^{\infty }p_{k} = 1 + K_{\alpha }\dfrac{\Delta t}{h^{\alpha }} \sum_{k=-\infty }^{\infty } w_{k} 
= 1 + K_{\alpha }\dfrac{\Delta t}{h^{\alpha }} \left[ w_{0} + \sum_{k=1 }^{\infty } (w_{k} + w_{-k}) \right] \text{.}
\label{e_FDM_suma}
\end{equation}
Substituting values from expressions (\ref{e_n01_w}) as well as~(\ref{e_n12_w}) into $w_{k}$ and making 
many calculations in both cases we finally obtain
\begin{equation}
\sum_{k=-\infty }^{\infty }p_{k} = 1 \text{.}
\end{equation}

In order to determine the stability of the explicit scheme~(\ref{e_FDMp2}) \cite{Blee01,Farl01}
the coefficient $p_{0}$ in~(\ref{e_FDM_p}) for $k=0$ should be positive 
(the other coefficients $p_{k}$ are non-negative for arbitrary values $\alpha $ and $\theta $ which is easy to prove)
\begin{equation}
p_{0} =1+K_{\alpha }\dfrac{\Delta t}{h^{\alpha }} w_{0} > 0\text{.}
\end{equation}
Thus we fixed the maximum length of the time step $\Delta t$, substituting into $w_{0}$ expressions~(\ref{e_n01_w})
and~(\ref{e_n12_w}) respectively, as
\begin{equation}
\Delta t<\dfrac{-h^{\alpha }}{K_{\alpha }w_{0}} = \dfrac{2h^{\alpha }}{K_{\alpha } \left( c_{L}+c_{R}\right) } 
\left\{ 
\begin{array}{ll}
\dfrac{\Gamma \left( 2-\alpha \right) }
{2^{1-\alpha }\lambda _{1}-3\lambda _{1}+2}
& \text{ for }0<\alpha < 1\text{,} \\ 
\dfrac{\Gamma \left( 3-\alpha \right) }
{2^{2-\alpha }\left( 2-\lambda _{2}\right) +4\lambda _{2}-6 } & \text{ for }1<\alpha \leq 2 \text{.}
\end{array}
\right.
\end{equation}

Moreover, the initial condition~(\ref{e_init}) is introduced directly to every grid
node at~the first time step $t=t^{0}$. This determines the initial values 
of function~$C$ as
\begin{equation}
C_{i}^{0}=c_{0}\left( x_{i}\right) \text{.}
\label{e_FDM_init}
\end{equation}

It is not easy to apply an implicit scheme in unbounded domains because 
it generates infinite dimensions for all matrices. Therefore, one
usually seeks improved difference equations within an explicit scheme.

\subsubsection{The boundary-initial value problem}
The numerical solution~(\ref{e_FDMp2}), which included the unbounded domain 
$-\infty <x<\infty $,  has no practical implementations in computer simulations.

\begin{figure}[h]
\begin{center}
 \includegraphics[width=0.7\textwidth]{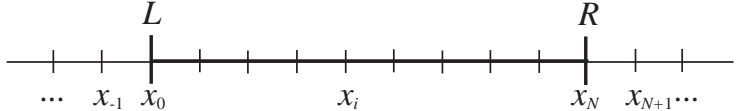}  
\end{center}
\caption{Grid nodes over space.}
\label{fig1}
\end{figure}

Here we try to solve this problem in the finite domain 
$\Omega :L\leq x\leq R$ with boundary conditions~(\ref{e_boundary}). We divide
domain $\Omega $ into $N$ sub-domains with the step $h = (R-L)/N$. Figure~\ref{fig1} shows
the modified spatial grid.
Here, we can observe additional 'virtual' points in the grid
located outside the lower and upper limits of the domain $\Omega $.
In order to introduce the Dirichlet boundary conditions, we propose
a numerical treatment which assumes the same values of function $C$
outside the domain limits as the values predicted on boundary nodes $x_{0}$ and $x_{N}$. 
\begin{equation}
C\left( x_{k},t\right) =\left\{ 
\begin{array}{ll}
C\left( x_{0},t\right) =g_{L}\left( t\right) & \text{ for }k<0\text{,} \vspace{0.1cm} \\ 
C\left( x_{N},t\right) =g_{R}\left( t\right) & \text{ for }k>N\text{.}
\end{array}
\right. 
\label{e_FDM_assume}
\end{equation}
Based on previous considerations we need to modify 
expressions~(\ref{e_n01_scheme}) and~(\ref{e_n12_scheme}) 
for the novel discretization of the Riesz-Feller derivative. Thus we have
\begin{equation}
_{x_{i}}D_{\theta }^{\alpha }C\left( x_{i},t\right)  
\approx \dfrac{1}{h^{\alpha }}
\left[ \sum\limits_{k=-i}^{N-i}C\left( x_{i+k},t\right) w_{k}
+ g_{L}\left( t\right) {s_{L}}_{i} 
+ g_{R}\left( t\right) {s_{R}}_{N-i} \right] ,
\label{e_n_mod}
\end{equation}
for $i=1,\ldots ,N-1$, where
\mathindent0cm
\begin{eqnarray}
{s_{L}}_{j} &=&{s_{L}}_{j}\left( \alpha ,\theta \right)
=\sum\limits_{k=-\infty }^{-j-1}w_{k}=\left\{ 
\begin{array}{ll}
c_{L}\cdot \overset{\sim }{r}_{j} & \text{ for }0<\alpha <1, \\ 
c_{L}\cdot \overset{\approx }{r}_{j} & \text{ for }1<\alpha \leq 2\text{,}%
\end{array}%
\right.  \label{e_sL} \\
{s_{R}}_{j} &=&{s_{R}}_{j}\left( \alpha ,\theta \right)
=\sum\limits_{k=j+1}^{\infty }w_{k}=\left\{ 
\begin{array}{ll}
c_{R}\cdot \overset{\sim }{r}_{j} & \text{ for }0<\alpha <1, \\ 
c_{R}\cdot \overset{\approx }{r}_{j} & \text{ for }1<\alpha \leq 2%
\end{array}%
\right.  \label{e_sP}
\end{eqnarray}
and
\begin{eqnarray}
\overset{\sim }{r}_{j} &=&\dfrac{\left( j+2\right) ^{1-\alpha }\lambda
_{1}+\left( j+1\right) ^{1-\alpha }\left( 2-2\lambda _{1}\right)
+j^{1-\alpha }\left( \lambda _{1}-2\right) }{2\Gamma \left( 2-\alpha \right) 
}, \\
\overset{\approx }{r}_{j} &=&\dfrac{\left( j+2\right) ^{2-\alpha }\left(
2-\lambda _{2}\right) +\left( j+1\right) ^{2-\alpha }\left( 3\lambda
_{2}-4\right) +j^{2-\alpha }\left( 2-3\lambda _{2}\right) +\left( j-1\right)
^{2-\alpha }\lambda _{2}}{2\Gamma \left( 3-\alpha \right) } \text{.} \notag \\
\end{eqnarray}

Modyfing expression~(\ref{e_FDMp1}) by substituting~(\ref{e_n_mod}) we obtain 
a finite difference scheme which is dependent on the weight factor $\sigma $. 
Here we assume
\begin{align}
g_{L}^{f+1/2} & = g_{L}\left( t^{f+1/2}\right) = g_{L}\left(
           \Delta t\left( f+\tfrac{1}{2}\right) \right), 
\label{e_gL}\\
g_{R}^{f+1/2} & = g_{R}\left( t^{f+1/2}\right) 
  = g_{R}\left( \Delta t\left( f+\tfrac{1}{2}\right) \right) 
\label{e_gR}
\end{align}
in order to simplify the numerical scheme. For internal nodes $x_{i}$, $i=1,\ldots ,N-1$
we have
\begin{align}
\dfrac{C_{i}^{f+1}-C_{i}^{f}}{\Delta t} = K_{\alpha }\dfrac{1}{h^{\alpha }}
 & \left[ \sum\limits_{k=-i}^{N-i}\left( \sigma C_{i+k}^{f}+\left( 1-\sigma
\right) C_{i+k}^{f+1}\right) w_{k} \right.   \notag
\\
&\left. +\ g_{L}^{f+1/2}{s_{L}}_{i}
+\ g_{R}^{f+1/2}{s_{R}}_{N-i} \right] 
\label{e_FDM_boundary}
\end{align}
and for boundary nodes $x_{0}$, $x_{N}$ we denote
\begin{align}
C_{0}^{f+1} &={\ g_{L}^{f+1/2},} \label{e_FDM_bL} \\
C_{N}^{f+1} &={\ g_{R}^{f+1/2}.} \label{e_FDM_bR}
\end{align}
The method is explicit for $\sigma =1$, partially implicit
for $0<\sigma <1$ and fully implicit for $\sigma =0$. 

The above scheme described by expressions~(\ref{e_FDM_boundary})-(\ref{e_FDM_bR}) 
can written in the matrix form as
\begin{equation}
\mathbf{A}\cdot \mathbf{C^{f+1}}=\mathbf{B}\text{,}
\label{e_FDM_matrix}
\end{equation}
where
\begin{equation}
\mathbf{A}=\left[ 
\begin{array}{cccccccc}
1 & 0 & 0 & 0 & \ldots  & 0 & 0 & 0 \\ 
a_{-1} & 1+a_{0} & a_{1} & a_{2} & \ldots  & a_{N-3} & a_{N-2} & a_{N-1} \\ 
a_{-2} & a_{-1} & 1+a_{0} & a_{1} & \ldots  & a_{N-4} & a_{N-3} & a_{N-2} \\ 
a_{-3} & a_{-2} & a_{-1} & 1+a_{0} & \ldots  & a_{N-5} & a_{N-4} & a_{N-2} \\ 
a_{-4} & a_{-3} & a_{-2} & a_{-1} & \ldots  & a_{N-6} & a_{N-3} & a_{N-4} \\ 
\vdots  & \vdots  & \vdots  & \vdots  & \ddots  & \vdots  & \vdots  & \vdots \\ 
a_{-N+2} & a_{-N+3} & a_{-N+4} & a_{-N+5} & \ldots  & 1+a_{0} & a_{1} & a_{2} \\ 
a_{-N+1} & a_{-N+2} & a_{-N+3} & a_{-N+4} & \ldots  & a_{-1} & 1+a_{0} & a_{1} \\ 
0 & 0 & 0 & 0 & \ldots  & 0 & 0 & 1
\end{array}
\right] , \quad
\mathbf{B}=\left[ 
\begin{array}{c}
{\ g_{L}^{f+1/2}} \\ 
b_{1} \\ 
b_{2} \\ 
b_{3} \\ 
b_{4} \\ 
\vdots  \\ 
b_{N-2} \\ 
b_{N-1} \\ 
{\ g_{R}^{f+1/2}}
\end{array}
\right] ,
\label{e_FDM_matrixes}
\end{equation}
with 
\begin{eqnarray}
a_{j} &=&\left(\sigma -1 \right) K_{\alpha }\dfrac{\Delta t}{h^{\alpha }}
w_{j} \quad \text{ for }j=-N+1,\ldots ,N-1\text{,}
\label{e_matrix_a} \\
b_{j} &=&C_{j}^{f}+K_{\alpha }\dfrac{\Delta t}{h^{\alpha }}
\left[ g_{L}^{f+\frac{1}{2}}{s_{L}}_{j}
+g_{R}^{f+\frac{1}{2}}{s_{R}}_{N-j}   
+ \sigma \sum\limits_{k=-j}^{N-j}C_{i+k}^{f}w_{k} \right]  
 \notag \\
&& \hspace{5cm} \text{for }j=1,\ldots ,N-1\text{.}
\label{e_matrix_b}
\end{eqnarray}
and $\mathbf{C^{f+1}}$ is the vector of unknown function values at time $t^{f+1}$.

A particular case of the scheme~(\ref{e_FDM_boundary}) is the explicit scheme for $\sigma =1$.
This case may be simplified to 
\begin{equation}
C_{i}^{f+1}=\left\{ {
\begin{array}{ll}
g_{L}^{f+1/2} & \hspace{-0cm}\text{for }i=0, \\ 
\displaystyle\sum\limits_{k=-i}^{n-i}C_{i+k}^{f}p_{k}
+ K_{\alpha }\dfrac{\Delta t}{h^{\alpha }}
\left( g_{L}^{f+1/2}{s_{L}}_{i}
      +g_{R}^{f+1/2}{s_{R}}_{N-i} \right)  & \text{for }i=1,\ldots ,N-1, \\ 
g_{R}^{f+1/2} & \hspace{-0cm}\text{for }i=N,
\end{array}
}\right. 
\label{e_FDM_explicit}
\end{equation}
where $p_{k}$ is defined by formula~(\ref{e_FDM_p}).
We can observe that boundary conditions influence on all the values
of the function at every node. Unlike the classical second derivative
which is approximated locally, the Riesz-Feller and other fractional derivatives 
accumulate all values of the function at the domain points.

The skewness parameter $\theta$ has a significant influence on the solution.
Assuming $\theta \rightarrow \pm 1^{\pm}$ and $\alpha \rightarrow 1^{\pm}$  one can obtain
the classical hyperbolic equation called the first order wave equation (the transport equation).
In this case our scheme tends to the known Euler's forward-time and backward-space/forward-space scheme.
For $\alpha = 2$, $\theta = 0$ our scheme is the same
as the forward-time and central-space scheme~\cite{Ames01,Hoff01}.


\section{ Simulation results and their analysis}
\begin{figure}[h]
\begin{center}
 \includegraphics[width=0.7\textwidth]{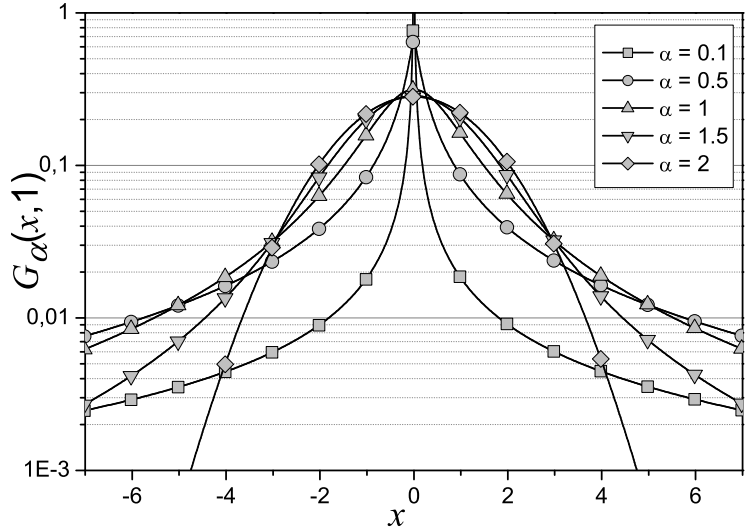}  
\end{center}
\caption{Probability density functions over space 
for different values of parameter $\alpha$.}
\label{fig2}
\end{figure}

In this section we present the results of calculation obtained by
our numerical approach. We try to simulate the evolutions of
the probability density function over time
for $\alpha \in \left\langle 0.1,0.5,1^{-},1.5,2\right\rangle $.
Assuming the domain $\left[ -10,10\right] $ we divided this domain into 1000 
subintervals ($h=0.02$). We assumed the boundary conditions
as $g_{L}\left( t\right) =g_{R}\left( t\right) =0$.
The initial condition $u(x,0^{+}) =\delta (x)$ is approximated by $C_{500}^{0}=1/h$ 
and $C_{i}^{0}=0$, for $i\neq 500$. 
Fig.~\ref{fig2} shows graphically the probability density function $G_{\alpha}(x,t)$
over space after $t=1$s in the limited interval $\left[ -7,7\right] $
(in the logarithmic scale).
It should be noted that for $\alpha =2$ our solution roughly estimates  the Gaussian 
probability density function: 
$G_2(x,1) = \left( 2\sqrt{\pi }\right) ^{-1}\exp \left( -x^{2}\right) $.
For $\alpha =1^{\pm }$ this solution becomes the Cauchy probability density function:
$G_1(x,1) = 1/\left( \pi \left( 1+x^{2}\right) \right) $.  

\begin{figure}[h]
\begin{center}
\includegraphics[width=0.7\textwidth]{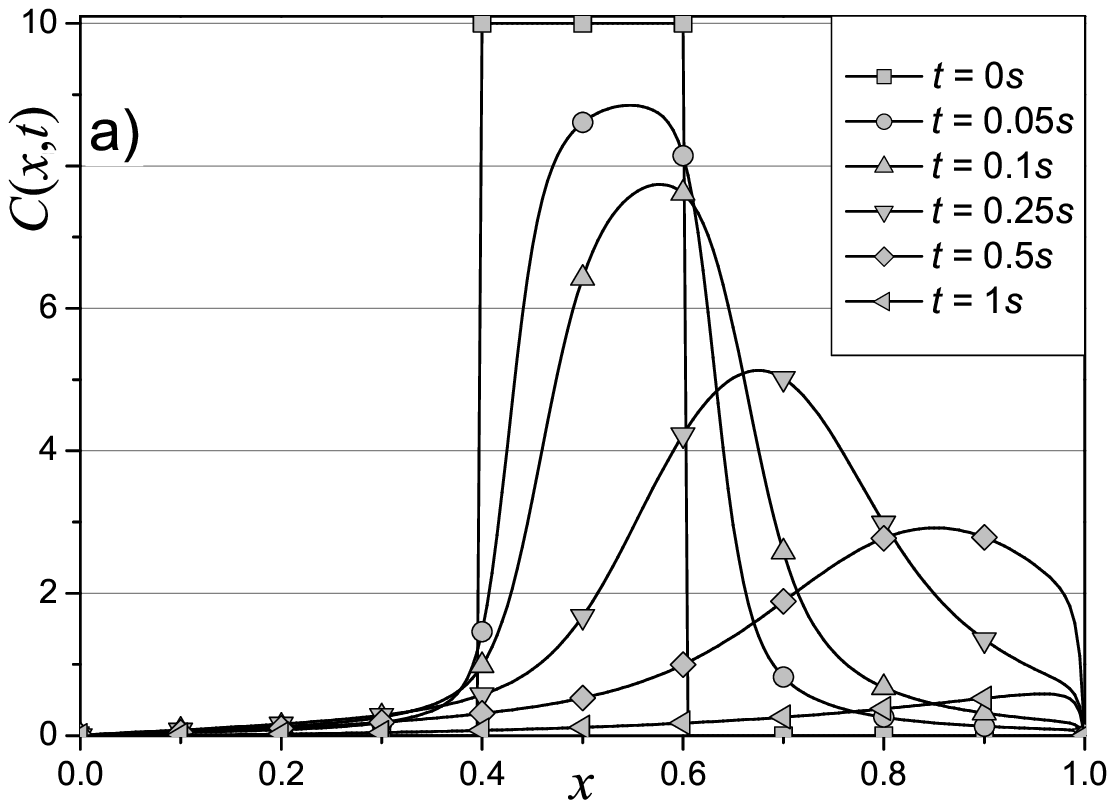}
\includegraphics[width=0.7\textwidth]{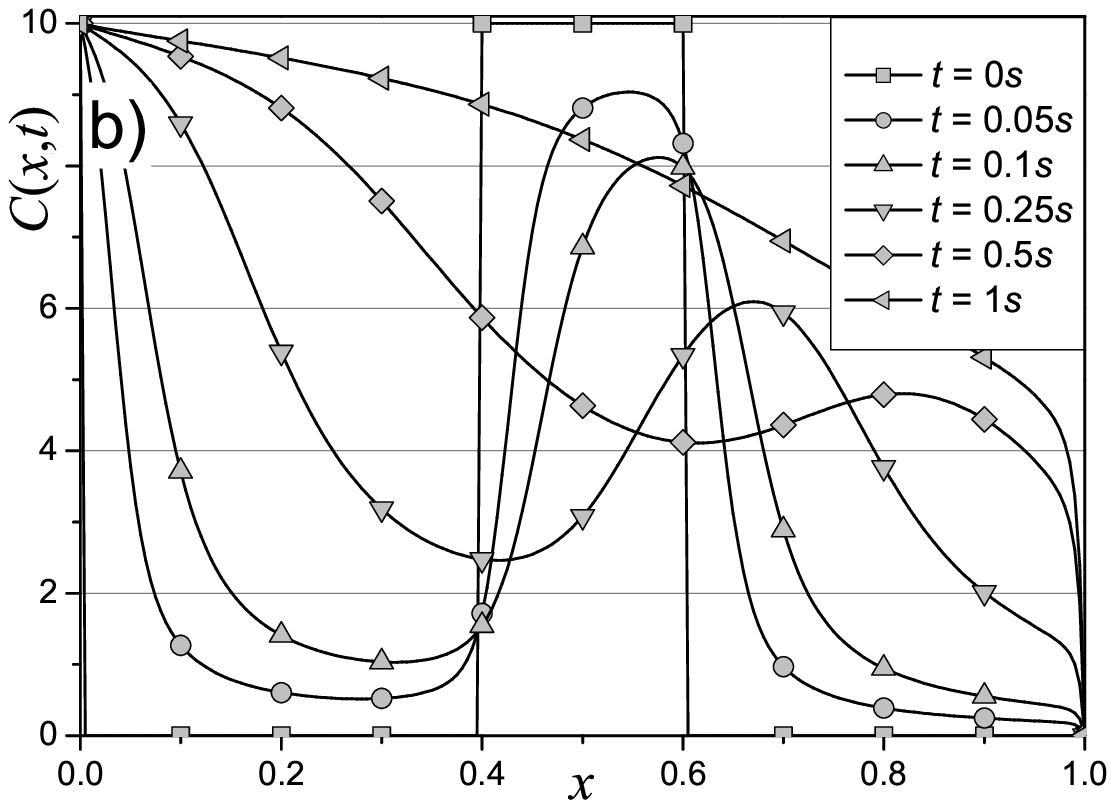}
\end{center}
\caption{Numerical solution to Eq.~(\ref{e_ad}) over space for $\alpha =0.9$,
$\theta = -0.7$ and for different values of boundary conditions.}
\label{fig3}
\end{figure}

\begin{figure}[!h]
\begin{center}
\includegraphics[width=0.7\textwidth]{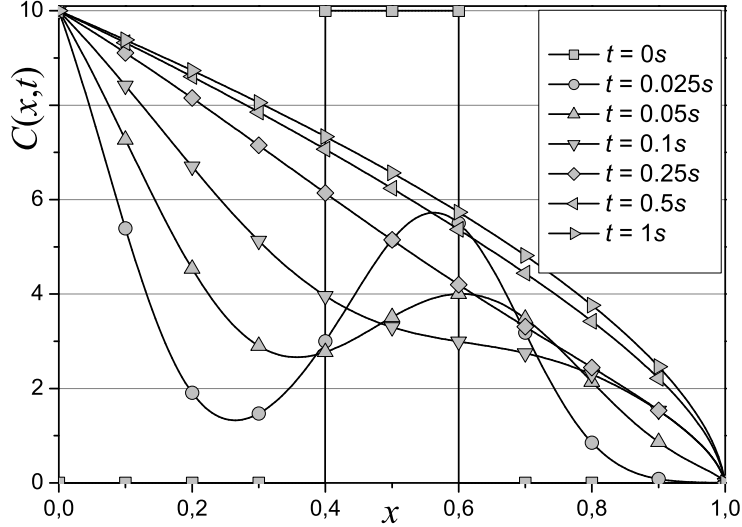}
\end{center}
\caption{Numerical solution to Eq.~(\ref{e_ad}) 
over space for $\alpha =1.6$, $\theta =-0.4$.}
\label{fig4}
\end{figure}

In the next two examples we show the influence of parameters $\alpha $ and 
$\theta $ on the solution. In both examples we consider the domain 
$\Omega :\left[ 0,1\right] $ with initial condition 
\begin{equation}
c_{0}\left( x\right) =\left\{ 
\begin{array}{ll}
10 & \text{for }x\in \left[ 0.4,0.6\right] , \\ 
0 & \text{for }x\in \left[ 0,0.4\right) \cup \left( 0.6,1\right] .
\end{array}
\right. 
\end{equation}%
Fig. \ref{fig3} shows the solution to function $C \left( x,t \right) $ for $\alpha =0.9$ and 
$\theta =-0.7$ with boundary conditions $g_{L}\left( t\right) =g_{R}\left( t\right) =0$ on Fig.~\ref{fig3}(a) 
and $g_{L}\left( t\right) = 10$, $g_{R}\left( t\right) =0$ on the Fig.~\ref{fig3}(b), 
at different moments of time.
Whereas Fig.~\ref{fig4} presents plots of $C\left( x,t\right) $, for $\alpha =1.6$, $%
\theta =-0.4$ and\ $g_{L}\left( t\right) =10$, $g_{R}\left( t\right) =0$.

\section{ Conclusions}
In summary we proposed the FFDM for fractional diffusion equations in which
the Riesz-Feller fractional derivative is included.
We analysed the anomalous diffusion equation in linear form in order 
to compare numerical results with the analytical solution. 
We obtained implicit and explicit FDM schemes 
which may generalise classical schemes of FDM. Moreover, our solution for $\alpha = 2$
is the same as the classical finite difference method. 
We hope that this numerical approach will be successfully applied to 
fractional partial differential equations having more complex forms, i.e. non-linear forms.

Analysing the graphs included in this work,
we observe that for the case $\alpha < 2$ (the L\'{e}vy flight)
diffusion is slower than classical diffusion (Brownian motion) 
at the initial time steps. Nevertheless, the probability density function 
generates a long tail of distribution in the long time limit. 
This can be associated with rare and extreme events which are
characterized by very large arbitrary values of particle jumps.

Analysing changes in the skewness parameter $\theta$ we observed interesting
behaviour in the numerical solution to Eq.~(\ref{e_ad}).
For $\alpha \rightarrow 1^{+}$ and $\theta \rightarrow \pm 1^{+}$
we obtained the first order wave equation. Assuming $\theta \in (0,1)$ 
(with restrictions to the order $\alpha$) we generated a class of non-symmetric 
probability density functions. 

The proposed numerical scheme creates a bridge between Gaussian and Cauchy processes.
Our scheme is also a~bridge between diffusion and transport phenomena.

\end{document}